\documentclass[
11pt]{amsart}

\usepackage{amsmath}
\usepackage{amssymb}
\usepackage{latexsym}
\usepackage{amsfonts}

\begin{document}

\newcommand{\ci}[1]{_{ {}_{\scriptstyle #1}}}
\newcommand{\ti}[1]{_{\scriptstyle \text{\rm #1}}}

\newcommand{\norm}[1]{\ensuremath{\|#1\|}}
\newcommand{\abs}[1]{\ensuremath{\vert#1\vert}}
\newcommand{\nm}{\,\rule[-.6ex]{.13em}{2.3ex}\,}

\newcommand{\p}{\ensuremath{\partial}}
\newcommand{\pr}{\mathcal{P}}

\newcounter{vremennyj}

\newcommand\cond[1]{\setcounter{vremennyj}{\theenumi}\setcounter{enumi}{#1}\labelenumi\setcounter{enumi}{\thevremennyj}}

\newcommand{\pbar}{\ensuremath{\bar{\partial}}}
\newcommand{\db}{\overline\partial}
\newcommand{\D}{\mathbb{D}}
\newcommand{\T}{\mathbb{T}}
\newcommand{\C}{\mathbb{C}}
\newcommand{\N}{\mathbb{N}}
\newcommand{\bP}{\mathbb{P}}

\newcommand\cE{\mathcal{E}}
\newcommand\cP{\mathcal{P}}
\newcommand\cC{\mathcal{C}}
\newcommand\cH{\mathcal{H}}
\newcommand{\cU}{\mathcal{U}}

\newcommand{\be}{\mathbf{e}}
\newcommand{\bx}{\mathbf{x}}

\newcommand{\la}{\lambda}

\newcommand{\td}{\widetilde\Delta}

\newcommand{\tto}{\!\!\to\!}
\newcommand{\wt}{\widetilde}
\newcommand{\shto}{\raisebox{.3ex}{$\scriptscriptstyle\rightarrow$}\!}

\newcommand{\La}{\langle }
\newcommand{\Ra}{\rangle }
\newcommand{\ran}{\operatorname{Ran}}
\newcommand{\rank}{\operatorname{Rank}}
\newcommand{\tr}{\operatorname{tr}}
\newcommand{\codim}{\operatorname{codim}}
\newcommand\clos{\operatorname{clos}}
\newcommand{\spn}{\operatorname{span}}
\newcommand{\essinf}{\operatorname{essinf}}
\newcommand{\vf}{\varphi}
\newcommand{\f}[2]{\ensuremath{\frac{#1}{#2}}}


\newcommand{\entrylabel}[1]{\mbox{#1}\hfill}

\newenvironment{entry}
{\begin{list}{X}%
  {\renewcommand{\makelabel}{\entrylabel}%
      \setlength{\labelwidth}{55pt}%
      \setlength{\leftmargin}{\labelwidth}
      \addtolength{\leftmargin}{\labelsep}%
   }%
}%
{\end{list}}



\numberwithin{equation}{section}

\newtheorem{thm}{Theorem}[section]
\newtheorem{lm}[thm]{Lemma}
\newtheorem{cor}[thm]{Corollary}
\newtheorem{prop}[thm]{Proposition}

\theoremstyle{remark}
\newtheorem{rem}[thm]{Remark}
\newtheorem*{rem*}{Remark}

\title[Analytic Projections and the Corona Problem]{Analytic projections, Corona Problem and geometry of holomorphic vector bundles}
\author{Sergei Treil}
\thanks{The work of S.~Treil was supported by the National Science Foundation under Grant  DMS-0501065}
\address{Sergei Treil, Department of Mathematics\\ Brown University\\ 151 Thayer Street Box 1917\\ Providence, RI USA  02912}
\email{treil@math.brown.edu}

\author{Brett D. Wick}

\address{Brett D. Wick, Department of Mathematics\\ University of South Carolina\\ LeConte College\\ 1523 Greene Street\\ Columbia, SC USA 29208}
\email{wick@math.sc.edu}

\begin{abstract}
The main result of the paper is a theorem giving a sufficient condition for the existence of a bounded analytic projection onto a holomorphic family of generally infinite dimensional subspaces (a holomorphic sub-bundle of a trivial bundle). This sufficient condition is also necessary in the case of finite dimension or codimension of the bundle. A simple lemma of N.~Nikolski connects the existence of a bounded analytic projection with the Operator Corona Problem (existence of a bounded analytic left inverse for an operator-valued function), so as corollaries of the main result we obtain new results about the Operator Corona Problem.  In particular, we find a new sufficient condition, a complete solution in the case of finite codimension, and a solution of the generalized Corona Problem. 
\end{abstract}

\maketitle
\setcounter{tocdepth}{1}
\tableofcontents
\section*{Notation}
\begin{entry}
\item[$:=$] equal by definition;\\
\item[$\C$] the complex plane;\\
\item[$\D$] the unit disk, $\D:=\{z\in\C:\abs{z}<1\}$;\\
\item[$\T$] the unit circle, $\T:=\p\D=\{z\in\C:\abs{z}=1\}$;\\
\item[$\mu$] measure on $\D$ with $d\mu:=\f{2}{\pi}\log\f{1}{\abs{z}}dxdy$;\\
\item[$dm$] normalized Lebesgue measure on $\T$;\\ 

\item[$\p, \pbar$] $\p$ and $\pbar$-operators, 
$\p :=\frac12(\frac{\p}{\p x} - i \frac{\p}{\p y})$, 
$\pbar :=\frac12(\frac{\p}{\p x} + i \frac{\p}{\p y})$;

\item[$\Delta$] ``normalized'' Laplacian, $\Delta :=\p \pbar=
\frac14(\frac{\p^2}{\p x^2} + \frac{\p^2}{\p y^2})$;

\item[$\norm{\cdot}, \nm\cdot \nm$]  norm; since we are dealing with matrix- and operator-valued functions we will use the symbol $\|\,\cdot\,\|$ (usually with a
subscript) for the norm in a functions space, while $\nm\,\cdot\,\nm$ is used for the
norm in the underlying vector (operator) space.
 Thus for a vector-valued function
$f$ the symbol
$\|f\|_2$ denotes its
$L^2$-norm, but the symbol $\nm f\nm$ stands for the scalar-valued function whose
value at a point $z$ is the norm of the vector (operator) $f(z)$;  \\

\item[$H^2$, $H^\infty$] Hardy classes of analytic functions, 
$$
H^p := \left\{ f\in L^p(\T) : \hat f (k) := \int_\T f(z) z^{-k}
\frac{|dz|}{2\pi} = 0\ \text{for } k<0\right\}. 
$$
Hardy classes can be identified with spaces of analytic functions in the unit disk
$\D$.  In particular, $H^\infty$ is the space of all bounded
analytic in $\D$ functions; \\

\item[$H^2_E$] vector-valued  Hardy class $H^2$ with values in $E$; 

\item[$L^\infty_{\!E_*\shto E}$] class of bounded functions on the unit
circle $\T$ whose values are bounded operators from $E_*$ to $E$;
\item[$H^\infty_{\!E_*\shto E}$] operator Hardy class  of bounded analytic
 functions whose values are bounded  operators from $E_*$ to $E$;
$$
\|F\|_\infty := \sup_{z\in \D}
\nm F(z)\nm=\underset{\xi\in
\T}{\operatorname{ess sup}}\nm F(\xi)\nm;
$$

\item[$H^{\infty}_{E_*\to E}(\Omega)$] operator Hardy class of bounded analytic functions on an arbitrary domain $\Omega$,    $\displaystyle\norm{F}_{\infty}:=\sup_{z\in\Omega}\nm F(z)\nm$. 
\end{entry}

Throughout the paper all Hilbert spaces are assumed to be separable. We
always assume that in any Hilbert space, an orthonormal basis is fixed, so that
any operator $A:E\to E_*$ can be identified with its matrix. Thus, besides
the usual involution $A\mapsto A^*$ ($A^*$ is the adjoint of $A$), we 
have two more: $A\mapsto A^T$ (transpose of the matrix) and $A\mapsto
\overline A$ (complex conjugation of the matrix), so $A^* =(\overline
A)^T =\overline{A^T}$. Although everything in the paper can be presented in
invariant, ``coordinate-free'' form, use of transposition and complex
conjugation makes the notation easier and more transparent. 

\setcounter{section}{-1}

\section{Introduction and Main Results}

The Operator Corona Problem is to find a (preferably local) necessary and
sufficient condition for a bounded operator-valued function $F\in
H^\infty_{E_*\to E}$ to have a left inverse in $H^\infty_{E_*\shto E}$, i.e.,~a function
$G\in H^\infty_{\!E\to E_*}$ such that 
\begin{equation}
\label{B}
\tag{B}
G(z)F(z)\equiv I \qquad \forall z\in \D .
\end{equation}
In the literature, such equations are sometimes called Bezout equations, and
``B'' here is for Bezout. The simplest necessary condition for \eqref{B} is 
\begin{equation}
\label{C}
\tag{C}
F^*(z)F(z) \ge \delta^2 I, \qquad \forall z\in \D\qquad (\delta>0)
\end{equation}
(the tag ``C'' is for Carleson).   If condition \eqref{C} implies \eqref{B}, we say that the Operator 
Corona Theorem holds. In the particular case when $F$ is a column  $F=(f_1, f_2, \ldots, f_n)^T$ the Operator Corona Theorem is just the classical Carleson Corona Theorem.

The Operator Corona Theorem plays an important role in different areas of
analysis; in particular, in Operator Theory (angles between invariant subspaces,
unconditionally convergent spectral decompositions, see \cite{Nik-book-v1,
Nik-shift, Treil-OTAA-89, Treil-subspaces}) as well as in Control Theory and
other applications.

\subsection{Motivations}

There are several motivations for this paper. The first one is that in the matrix case, all the information about the Corona Problem is encoded in the analytic family of subspaces (a holomorphic vector bundle) $ \ran F(z)$, $z\in \D$. 

Let us explain this. Let $F\in H^\infty_{E_*\shto E}$ be the Corona Data satisfying $F^*F\ge\delta^2 I$ and let $F= F\ti i F\ti o$ be its inner-outer factorization. Then the outer part $F\ti o$ of $F$ is invertible and $\| F\ti o^{-1}\|_\infty\le 1/\delta$. Therefore, the Corona Problem for $F$ is equivalent to the Corona Problem for its inner part $F\ti i$. But, in this case the inner part $F\ti i$ of $F$ can be restored completely from the analytic family of subspaces $ \ran F(z)$, $z\in \D$.  Namely, one has to consider the $z$-invariant subspace $\cE := \{f\in H^2_E: f(z)\in \ran F(z) \ \forall z\in \D\}$, and the corresponding inner function is exactly $F\ti i$ (to see that $\cE \subset F\ti i H^2_{E_*}$ we can use the fact that $F\ti i^* F\ti i\ge\delta_1^2 I$; the opposite inclusion is trivial).  

So, a natural question arises: Is it possible to characterize condition \eqref{C} (or \eqref{B}) in purely geometric terms, i.e.,~in terms of the family of subspaces $\ran F(z)$, $z\in \D$? In this paper, such a characterization is given. 

Another motivation is the so-called \emph{Codimension One Conjecture}. It is trivial in the ``square case'', when the operators $F(z)$ are onto for all $z\in \D$, that condition \eqref{C} implies $F$ is invertible in $H^\infty$.  In this special case the Corona Problem is trivial. 

But, what happens if we consider the ``almost square'' case, i.e., the case when $\codim F(z)=1$ for all $z\in\D$? It looked plausible and it was conjectured by N.~Nikolski and the first author, that in this case the Operator Corona Theorem holds.  Besides the na\"\i ve reason that this case is close to the ``square'' one, there are additional reasons for this conjecture. One is that for any $n\times k$ Corona data $F$, one can canonically associate another set Corona data $\wt F$ of size $n\times (n-k)$ and the best possible norms of the $H^\infty$  left inverses of $F$ and $\wt F$ coincide.

It turns out that the Codimension One Conjecture failed, and failed spectacularly! An infinite dimensional counterexample was constructed in \cite{Tre-rank-one}, as well as lower bounds for the solution in the finite dimensional ``codimension one'' case. These lower bounds blow up as the dimension grows and they are almost optimal in that they are very close to the upper bound for the solution obtained recently by T.~Trent \cite{Trent}, see also \cite{TreilWick}. This means that from the point of view of the Operator Corona Theorem, the codimension one case is as bad (or almost as bad) as the general one.  

However, the codimension one case looks significantly simpler than the general one, and despite the failure of the Codimension One Conjecture, one can still hope to obtain a simple necessary and sufficient condition for the solvability of the Corona Problem. 

In this paper, we obtain a simple necessary and sufficient condition for the left invertibility of an operator-valued function in the case of finite codimension, not just codimension one.  

Finally, the following surprising lemma of N.~Nikolski was the last motivation for this paper.  Let $\Omega $ be a domain in $\C^n$, $n\ge1$ (in fact we can let $\Omega$ be a manifold). 
\begin{lm}[Nikolski's Lemma]
\label{l0.1}
Let $F\in H^\infty_{E_*\shto E}(\Omega)$ satisfy
$$
F^*(z)F(z)\ge \delta^2 I, \qquad \forall z\in \Omega. 
$$
Then $F$ is left invertible in $H^\infty_{E_*\shto E}(\Omega)$ (i.e.,~there exists $G\in
H^\infty_{E\shto E_*}(\Omega)$ such that $GF\equiv I$) if and only if there exists a
function $\cP\in H^\infty_{E\shto E}(\Omega)$ whose values are projections (not
necessarily orthogonal) onto
$F(z)E$  for all $z\in\Omega$.

Moreover, if such an analytic projection $\cP$ exists, one can find a left inverse $G\in H^\infty_{E\shto E_*}(\Omega)$ satisfying $\|G\|_\infty\le \delta^{-1} \|\cP\|_\infty$.
\end{lm}

This Lemma connects the two problems discussed above and plays an important part in their solution. 

\subsection{Main Results}

Instead of considering families of subspaces, we  consider more ``analytic'' objects; namely, the families of orthogonal projections $\Pi(z)$ onto these subspaces. The function $\Pi(z)$ is not analytic, except in the trivial case of a constant function.  The fact the family of subspaces $\ran \Pi(z)$ is an analytic family (a holomorphic vector bundle) is expressed by the identity $\Pi \p \Pi =0$. 

Let us now list the main results of the paper. Recall that $\Delta$ is the ``normalized'' Laplacian, $\Delta := \p\pbar = \frac14 \left(\frac{\p^2}{\p x^2} + \frac{\p^2}{\p y^2}\right)$. 

\begin{thm}[Main Result]
\label{t0.2}
Let $\Pi:\D \to B(E)$ be a $\cC^2$ function whose values are orthogonal projections in $E$ satisfying $\Pi\p\Pi =0$. Assume that there exists a bounded non-negative subharmonic function $\vf$ such that 
\begin{equation}
\label{0.1}
\Delta \vf (z) \ge \nm\p\Pi (z)\nm^2\qquad \forall z\in \D. 
\end{equation}
Then there exists a bounded analytic projection onto $\ran \Pi(z)$, i.e.,~a function $\cP\in H^\infty_{E\shto E}$ such that $\cP(z)$ is a projection onto $\ran\Pi(z)$ for all $z\in \D$.  

Moreover, if $0\le \vf(z)\le K$ for all $z\in\D$, then one can find $\cP$ satisfying 
$$
\|\cP\|_\infty \le 1+ 2  \sqrt{(Ke^{K+1}+1) Ke^{K+1} }.
$$
\end{thm}
By $\cC^2$ we mean twice continuously differentiable in the operator norm topology. This assumption can probably be relaxed, but in this paper we want to avoid unnecessary technical details.

\begin{rem}
Condition \eqref{0.1} simply means the Green potential 
$$
\mathcal G(\la) := \frac2\pi\iint_\D \ln \left| \frac{z-\la}{1-\overline\la z} \right| \nm{\p\Pi(z)}\nm^2 dxdy
$$
is uniformly bounded in the disk $\D$. Integrating separately over a small hyperbolic neighborhood of $\la$ and the rest of $\D$, it is not hard to see that the uniform boundedness of $\mathcal G$ follows from the two conditions below:
\begin{align}
	&\nm{\p\Pi(z)}\nm\le C/(1-|z|); \hfill\hfill\hfill 
	\label{0.2}\\
  &\text{the measure } \nm{\p\Pi(z)}\nm^2 (1-|z|) dxdy \text{ is a Carleson measure. \hfill}   \label{0.3}
\end{align}
\end{rem}

\begin{rem}
\label{rem0.4}
Since $\Pi = \Pi^*$, $\pbar \Pi=(\p \Pi)^*$, in conditions \eqref{0.1}, \eqref{0.2}, and \eqref{0.3} $\p$ can be replaced by $\pbar$. 
\end{rem}

The following proposition shows that in the case of finite dimension, or finite codimension, the above conditions  \eqref{0.2} and \eqref{0.3} (and therefore condition \eqref{0.1} of Theorem \ref{t0.2}) are necessary. Therefore, in this case \eqref{0.1} is equivalent to conditions \eqref{0.2} and  \eqref{0.3}.

\begin{prop}
\label{p0.5}
Suppose there exists a bounded analytic projection $\cP(z)$ onto $\ran\Pi(z)$, $z\in\D$. Assume also that 
either $\dim \ran\Pi(z)<\infty$ or \linebreak $\codim \ran\Pi(z)<\infty$%
.  Then conditions \eqref{0.2} and \eqref{0.3} (and therefore condition \eqref{0.1} of Theorem \ref{t0.2})  hold. 
\end{prop}

We note that, if there is a bounded analytic projection, both $\dim\ran\Pi(z)$ and  $\codim\ran\Pi(z)$ are constant for all $z\in \D$.  One of the main corollaries of the above results is the following theorem.
\begin{thm}[Operator Corona Theorem]
\label{t0.4}
Let $F\in H^\infty_{E_*\shto E}$ satisfy the Corona Condition $F^*F\ge\delta^2 I$. Assume also that the orthogonal projections $\Pi(z)$ onto $\ran F(z)$ satisfy assumption \eqref{0.1} of Theorem \ref{t0.2} (or conditions \eqref{0.2} and \eqref{0.3}). 
Then $F$ has a holomorphic left inverse $G\in H^\infty_{E\shto E_*}$. 

Moreover, if the function $\vf$ from condition \eqref{0.1} satisfies 
$$
0\le \vf(z)\le K\qquad \forall z\in\D, 
$$
then one can find the left inverse $G$ satisfying 
$$
\|G\|_\infty \le \tilde\delta^{-1}\left( 1+ 2  \sqrt{(Ke^{K+1}+1) Ke^{K+1} }\right),
$$
 where $\tilde\delta:=\essinf\{\nm F(z)e\nm : z\in \T, e\in E_*, \nm e\nm=1\}$.  
\end{thm}

Note that $\delta\le \inf\{\nm F(z)e\nm : z\in \D, e\in E_*, \nm e\nm=1\}$, so considering boundary values of $F(z)e$, one can easily see that $\tilde\delta\ge \delta$. 

Theorem \ref{t0.4} with $\delta$ instead of $\tilde\delta$ in the estimate of $\|G\|_\infty$ is an immediate corollary of Lemma \ref{l0.1} and Theorem \ref{t0.2}. The estimate with $\tilde\delta$ will be explained later in Section \ref{s1}, see Remark \ref{r1.1} there.

If $F\in H^\infty_{E_*\shto E}$ satisfies $F^*F\ge \delta^2 I$ and $\dim E_*<\infty$, an easy computation (see Proposition \ref{p3.1} below) shows that the orthogonal projection $\Pi(z)$ satisfies conditions \eqref{0.2} and \eqref{0.3} (and thus condition \eqref{0.1}). Therefore, the  Operator Corona  Theorem in the case $\dim E_*<\infty$ follows immediately from Theorem \ref{t0.4}. However, this result has been known for a long time as the Fuhrmann--Vasyunin Theorem, see \cite{Tolo}; see also \cite{Trent} or \cite{TreilWick} for the modern treatment with better estimates. 

Probably the most important new and non-trivial corollary is the following theorem, solving the operator Corona Problem in the case of finite codimension. 

\begin{thm}[Finite Codimension Operator Corona Problem]
\label{t0.6}
Let $F\in H^\infty_{E_*\shto E}$ satisfy the Corona Condition $F^*F\ge\delta^2 I$, and let $\codim\ran F(z)<\infty$. Then $F$ has a bounded analytic left inverse if and only if the orthogonal projections $\Pi(z)$ onto $\ran F(z)$ satisfy assumption \eqref{0.1} of Theorem \ref{t0.2} (or conditions \eqref{0.2} and \eqref{0.3}). 
\end{thm}

Finally, one can use Theorem \ref{t0.2} to obtain results about generalized inverses of $F$ in the case when $\ker F(z)\ne \{0\}$. Let us recall (cf \cite[Chapter 12]{Lancaster-ToM-1985} for the matrix case) that a generalized  inverse of an operator $A$ is an operator $B$ such that 
$$
ABA  = A \textnormal{ and } BAB=B.  
$$
The condition $ABA=A$ simply means that the linear operator $B$, defined on the whole space, gives a solution $x=Bb$ of the equation $Ax =b$ whenever such a solution exists. Note that if $A$ is left (right) invertible a generalized inverse $B$ is simply a left (right) inverse.

A generalized inverse is usually  not unique. 
However, imposing the additional requirement that the operators $AB$ and $BA$ are self-adjoint, one can make the generalized inverse unique.  In this case the unique generalized inverse is called the \emph{Moore--Penrose} inverse.  

Often in the literature, the term ``generalized inverse'' is  used for the Moore--Penrose inverse, cf \cite{Campbell-GenInv-1991}.  We  want  to emphasize that  we are not dealing with the Moore--Penrose inverse here. Theorem \ref{gencorona} and Lemma \ref{GNL} below are trivially false for the Moore--Penrose inverse. 

While every non-zero matrix has a generalized inverse, this is not the case for the operators acting from one Hilbert space to another. It is well known and is easy to see that in this case, an operator $A$ has a generalized inverse if and only if 
$$
\nm A\bx \nm \ge\delta\nm \bx \nm \qquad \forall\bx \in (\ker A)^\perp
$$
for some $\delta>0$. 

Given $F\in H^\infty_{E_*\shto E}$, we are interested in finding a bounded analytic generalized inverse $G$, i.e., $G\in H^\infty_{E\shto E_*}$ that satisfies 
\begin{equation}
\label{B1}\tag{B1}
F(z)G(z)F(z)  =F(z), \quad G(z)F(z)G(z)=G(z) \qquad \forall z\in \D. 
\end{equation}
Clearly the condition
\begin{equation}
\label{C1}
\tag{C1}
\nm{F(z)\be}\nm\ge \delta \nm{ \be}\nm\qquad   \forall \be \in (\ker F(z))^\perp \quad \forall z\in \D ,
\end{equation}
is necessary  for the existence of a bounded  generalized inverse. 
Note, that if $\ker F(z)=\{0\}$, then condition \eqref{C1} is just the Carleson Condition \eqref{C}.

\begin{thm}[Generalized Operator Corona Theorem]
\label{gencorona}
Let $F\in H^\infty_{E_*\shto E}$ satisfy \eqref{C1}. Assume that the orthogonal projections $\Pi_R(z)$ onto $\ran F(z)$ and $\Pi_K(z)$ onto $\ker F(z)$ satisfy assumption \eqref{0.1} of Theorem \ref{t0.2} (or conditions \eqref{0.2} and \eqref{0.3}).  Then $F$ has a bounded analytic generalized inverse, i.e., there exists  $G\in H^\infty_{E\shto E_*}$ satisfying \eqref{B1}.

Moreover, if $\rank F(z)<\infty$ for all $z\in \D$ then conditions \eqref{0.2} and \eqref{0.3} (and therefore condition \eqref{0.1} of Theorem \ref{t0.2}) are necessary for the existence of a bounded analytic generalized inverse. 

\end{thm}

The main step in the proof of Theorem \ref{gencorona} is the following lemma, which holds for arbitrary domains in $\C^n$, $n\ge 1$ (and even for manifolds), not just for $\D$.

\begin{lm}[Generalized Nikolski Lemma]
\label{GNL}
Let $F\in H^\infty_{E_*\shto E}(\Omega)$ satisfy
$$
\nm F(z) x\nm\ge \delta \nm x\nm \quad\forall x\in(\ker F(z))^\perp \quad\forall z\in \Omega.
$$
The following statements are equivalent: 
\begin{enumerate}
\item $F$ has a generalized inverse in $H^\infty_{E_*\shto E}(\Omega)$; 

\item There exists $G\in
H^\infty_{E\shto E_*}(\Omega)$ such that $FGF\equiv F$ in $\Omega$;


\item There exist bounded analytic projections onto $\ran F(z)$ and $\ker F(z)$, i.e. there exists 
functions $\cP_R\in H^\infty_{E\shto E}(\Omega)$ and $\cP_K\in  H^\infty_{E_*\shto E_*}(\Omega)$ 
whose values are projections (not necessarily orthogonal) such that 
$$
\ran \cP_R (z) = \ran F(z) , \quad \ran \cP_K(z) =\ker F(z) \quad \forall z\in \Omega. 
$$
\end{enumerate}
\end{lm}
\begin{rem*}
The indices $R$ and $K$ in $\cP_R$ and $\cP_K$ stand for \emph{range} and \emph{kernel}. 
\end{rem*}

Sufficiency in Theorem \ref{gencorona} follows immediately from the above lemma and Theorem \ref{t0.2}. Necessity follows from  Proposition \ref{p3.1} below in Section \ref{s3}. 

\section{Nikolski's Lemma and Variations}
\label{s1}

\subsection{Proof of Nikolski's Lemma} The proof of Nikolski's Lemma (Lemma \ref{l0.1}) presented below is due to N.~Nikolski (personal communication), see also \cite{Tre-oper-corona}. The proof   is so simple that it is almost as surprising as the result!

\begin{proof}[Proof of Lemma \ref{l0.1}]
Let $F$ be left invertible in $H^\infty_{E_*\shto E}(\Omega)$ and let $G$ be one of its left inverses.
Define $\cP\in H^\infty_{E\shto E}(\Omega)$ by
$$
\cP(z) = F(z) G(z).
$$
Direct computation shows that $\cP^2 =\cP$, so the values of $\cP$ are projections.

Since $GF\equiv I$, 
$$
G(z)E = E_* \qquad \forall z\in \Omega.
$$ 
Therefore
$$
\cP(z) E = F(z)G(z) E = F(z) E_* \qquad \forall
z\in \Omega, 
$$
i.e.,~$\cP(z)$ is indeed a projection onto $F(z)E$. 

Suppose now that there exists a projection-valued function $\cP\in
H^\infty_{E\shto E}(\Omega)$, whose values are projections onto $F(z)E$ for all $z\in
\Omega$. We want to show that $F$ is left invertible in $H^\infty_{E_*\shto E}(\Omega)$.

First of all, let us notice that locally, in a neighborhood of each point $z_0\in\Omega$, the function $F(z)$ has an analytic left inverse. Indeed, if an operator $G_0:E_*\to E$ is a constant left inverse to the operator $F(z_0)$, i.e., if $G_0 F(z_0)=I$, then 
$$
G_0 F(z) = I - G_0 (F(z_0) -F(z) ). 
$$
So the inverse of $G_0F(z)$ is given by the analytic function 
$$
A(z) := \sum_{k=0}^\infty \left[ G_0 \cdot (F(z_0) -F(z) ) \right]^k
$$
defined in a neighborhood of $z_0$. So, $A(z)G_0$ is a local analytic left inverse of $F(z)$.

Fix  an arbitrary bounded left inverse $F^\dagger\in L^\infty_{E\shto E_*}(\Omega)$ of $F$, for example $F^\dagger = (F^*F)^{-1}F^*$. Note that we do not claim  here that $F^\dagger $ is analytic, only that it is bounded.

Note that for any other  left inverse $\wt G(z)$ of $F(z)$
\begin{equation}
\label{1.1}
\wt G(z) \bigm| \ran F(z)  = F^{\dagger}(z) \bigm| \ran F(z) 
\end{equation}
(action of a left inverse on  $\ran F(z)$ is uniquely defined).

Define  a function $G$ by 
$$
G(z) := F^\dagger(z)\cP(z), 
$$
which is well defined. It is easy to see that $G$ is bounded (since both $F^\dagger$ and $\cP$ are bounded) and that $G(z)F(z) \equiv I$. 

Let us show that $G$ is analytic. 
Fix a point $z_0\in\Omega$ and let $G_{z_0}(z)$ be a {\em local} analytic left inverse of $F(z)$ defined in a neighborhood of $z_0$. It follows from \eqref{1.1} that 
$$
G(z) := F^\dagger(z)\cP(z) = G_{z_0}(z)\cP(z) 
$$
in a neighborhood of $z_0$, so $G(z)$ is analytic there. Since $z_0$ is arbitrary, $G$ is analytic in $\Omega$. 

Finally, for $F^\dagger = (F^*F)^{-1}F^*$ we have $\|F^\dagger\|_\infty\le 1/\delta$, so $G=F^\dagger \cP$ satisfies  $\|G\|_\infty\le \delta^{-1}\|\cP\|_\infty$. 
\end{proof}

\begin{rem}
\label{r1.1}
For the case $\Omega=\D$ (or, more generally, when $\Omega$ is a domain such that the norm of an $H^\infty$ function is defined by the norm of its boundary values) the estimates of the norm of the left inverse $G$ can be improved to $\|G\|_\infty\le \tilde\delta^{-1}\|\cP\|_\infty$, where  $\tilde\delta:=\essinf\{\nm F(z)e\nm : z\in \T, e\in E_*, \nm e\nm=1\}$. 

Indeed, for $G=F^\dagger \cP$, where $F^\dagger = (F^*F)^{-1}F^*$, the equality 
$$
\lim_{z\to\zeta} G(z) = \lim_{z\to\zeta} F^\dagger(z) \lim_{z\to\zeta} \cP(z)
$$
holds for almost all $\zeta\in\T$ (here $\displaystyle\lim_{z\to\zeta}$ is the non-tangential limit in the strong operator topology). 

The $L^\infty(\T)$ norm of the boundary values of $F^\dagger$ can be estimated above by $\tilde\delta^{-1}$, and the norm of an $H^\infty$ function  equals to the norm of its boundary values. Therefore $\|G\|_\infty\le \tilde\delta^{-1}\|\cP\|_\infty$. 
\end{rem}

\subsection{Co-Outer Functions and Nikolski's Lemma on the Boundary }
\label{s1.2}

In Nikolski's Lemma we have the analytic projection $\cP(z)$ onto $\ran F(z)$ for all $z\in \D$. In the matrix case, $\dim E_*<\infty$, standard reasoning allows one to replace this condition by the existence of an $H^\infty$ projection onto $\ran F(z)$ a.e.~on $\T$. However, in the general case this is not true. In \cite{Tre-oper-corona} an example of a function $F\in H^\infty_{E_*\shto E}$ satisfying $F^*F\ge\delta^2 I$ on $\D$ but not left invertible in $H^\infty_{E_*\shto E}$ and such that $\ran F(z)=E$ a.e.~on $\T$ was constructed. However, there exists a bounded analytic projection (the trivial one, $\cP(z)\equiv I$) onto $\ran F(z)$ a.e.~on $\T$. By adding extra dimensions to $E$, one can make this projection non-trivial (different from $I$).

However, in the general case it is still possible to obtain a ``boundary'' version of Nikolski's Lemma. To do that, we need the notion of a co-outer function.

\subsubsection{Outer and Co-Outer Functions}

Let us recall that an operator-valued function $F\in H^\infty_{E_*\shto E}$ is called \emph{inner} if $F(z)$ is an isometry a.e.~on~$\T$, and $F$ is called \emph{outer} if $F H^2_{E_*}$ is dense in $H^2_E$ (or equivalently, if the set $\{Fp: p \text{ is a polynomial in }H^2_{E_*}\}$ is dense in $H^2_E$).   A function $F\in H^\infty_{E_*\shto E}$ is called \emph{co-outer} if the function $F^T$ is outer. 

Recall that we are working with fixed orthonormal bases in $E$ and $E_*$, so the transpose and complex conjugate of an operator, i.e.,~of its matrix, are well defined.

The following proposition gives several equivalent definitions of a co-outer function. 

\begin{prop}
\label{p1.2}
Let $F\in H^\infty_{E_*\to E}$. The following statements are equivalent:
\begin{enumerate}
	\item $F^T$ is outer; 
	\item \label{test2}$F^* ((H^2_{E})^\perp ) $ is dense in $(H^2_{E_*})^\perp$;
	\item $F^\sharp$ defined by $F^\sharp(z):= F^*(\overline z)$ is outer. 
\end{enumerate}
\end{prop}
\begin{proof}
By  the definition $F^T$ is outer if and only if $F^TH^2_{E}$ is dense in $H^2_{E_*}$, or, equivalently $F^T zH^2_{E}$ is dense in $zH^2_{E_*}$. Taking complex conjugates (we are assuming that orthonormal bases in $E$ and $E_*$ are fixed, so vectors can be identified with columns, and operators with  matrices), we see that the latter is equivalent to the fact that  $F^* ((H^2_{E})^\perp )$ is dense in $(H^2_{E_*})^\perp$, so \cond 1$\Longleftrightarrow$\cond 2. 

On the other hand, the transformation $\tau$, $(\tau f)(z) := f(\bar z)$ unitarily maps  $(H^2)^\perp$ onto $zH^2$, and $\tau (F^*) = F^\sharp$. 
Therefore, applying $\tau$ to property (2) we get that (2) holds if and only if $F^\sharp zH^2_E $ is dense in $zH^2_{E_*}$, i.e., that 
\cond 2$\Longleftrightarrow$\cond 3.
\end{proof}

It is an easy exercise to check, that for a co-outer function $F$, the operators $F(z)$ have trivial kernel for all $z\in \D$. 

To see that, let us recall that for an operator-valued function $G\in L^\infty_{E\shto E_*}$ the T\"oplitz operator $T_G:H^2_E \to H^2_{E_*}$ is defined by 
$$
T_G f = P_+(Gf) , \qquad f\in H^2_E, 
$$
where $P_+$ is the orthogonal projection (in $L^2(\T)$) onto $H^2$. 

Let $G=F^T\in H^\infty_{E\shto E_*}$. The function  $F$ is co-outer if and only if the T\"oplitz operator $T_G$ has dense range, or equivalently, if and only if $(T_G)^* = T_{G^*}$ has trivial kernel. 

Let $k_\la$ be the reproducing kernel of $H^2$, $k_\la(z) = 1/(1-\overline \la z)$. It is well known and easy to show that for $G\in H^\infty_{E\shto E_*}$
$$
T_{G^*} k_\la e = k_\la G^*(\la) e, \qquad \forall e\in E, \quad \forall \la\in \D. 
$$

Therefore, if $G$ is outer (i.e.~if $F$ is co-outer), $\ker G^*(\la)$ is trivial for all $\la\in \D$. 
But $G^*=\overline F$, so $\ker F(\la)$ is also trivial for all $\la\in \D$.

\medskip

The following proposition characterizes co-outer functions.

\begin{lm}
\label{l-1.3}
Let $F\in H^\infty_{E_*\shto E}$ satisfy 
$$
F^*(z)F(z) \ge\tilde \delta^2 I, \qquad \text{a.e.~on }\T
$$
for some $\tilde\delta>0$ (i.e.,~$F$ is left invertible in $L^\infty_{E_*\shto E}$).   

The function $F$ is co-outer if and only if the subspace 
$$
\cE =\cE_F:= \{ f\in H^2_E : f(z) \in \ran F(z) \ \text{a.e.~on } \T \}
$$
coincides with $FH^2_{E_*}$. 
\end{lm}

\begin{proof}

Let $F$ be co-outer. Since trivially  $FH^2_{E_*}\subset \cE$, one only needs to check the inclusion $\cE\subset FH^2_{E_*}$.  

Take $f\in \cE$. Since $F$ is left invertible in $L^\infty_{E_*\shto E}$, there exists a unique function $g\in L^2_{E_*}$ such that 
$$
f = F g. 
$$
Since $Fg=f \in\cE\subset  H^2_E$, we have that 
$$
0 = \La Fg, h\Ra = \La g, F^* h\Ra\qquad \forall h\in (H^2_E)^\perp. 
$$
Since $F$ is co-outer, the set $F^*(H^2_E)^\perp$ is dense in $(H^2_{E_*})^\perp$, so $g\in H^2_{E_*}$. But that means $f=Fg\in FH^2_{E^*}$.

To prove the opposite implication assume that $\cE = FH^2_{E_*}$. Assume that $F$ is not co-outer, so there exists a function $g\in (H^2_{E_*})^\perp \ominus F^* (H^2_{E})^\perp$, so  
$$
0=\La g, F^*h\Ra= \La Fg, h\Ra  \qquad \forall h \in (H^2_{E})^\perp. 
$$
Therefore $f:=Fg\in H^2_E$, and so $f\in \cE$. But, we assumed that $g\in (H^2_{E_*})^\perp$, which contradicts the assumption $\cE = FH^2_{E_*}$ (since $F$ is left invertible in $L^\infty_{E_*\shto E}$, the function $g$ is the unique $L^2_{E_*}$ solution of the equation $Fg=f$). 
\end{proof}

As mentioned at the beginning of Section \ref{s1.2}, the condition $F^*(z)F(z)\ge \delta^2 I$ for all $z\in \D$ does not generally mean that the function $F$ is co-outer. However, the stronger condition of left invertibility in $H^\infty$ implies that $F$ is co-outer. 
\begin{lm}
Any operator-valued function $F\in H^\infty_{E_*\shto E}$ which is left invertible in $H^\infty_{E_*\shto E}$ (i.e., there exists a $G\in H^\infty_{E\shto E_*}$ such that $GF\equiv I$ for all $z\in\D$) is co-outer.
\end{lm}
\begin{proof}
If $GF\equiv I$ then 
$$
(H^2_{E_*})^\perp = F^*G^*((H^2_{E_*})^\perp) \subset F^*((H^2_{E})^\perp) \subset (H^2_{E_*})^\perp, 
$$
so all the inclusions must be equalities. In particular, $F^*((H^2_{E})^\perp) = (H^2_{E_*})^\perp$, so by Proposition \ref{p1.2}, $F$ is co-outer. 
\end{proof}

\subsubsection{Nikolski's Lemma on the boundary.}
\begin{lm}
Let $F\in H^\infty_{E_*\shto E}$ be a co-outer function satisfying
$$
F^*(z)F(z) \ge \delta^2 I \qquad \text{a.e.~on }\T.
$$
Then $F$ is left invertible in $H^\infty_{E_*\shto E}$ if and only if there exist $\cP\in H^\infty_{E\shto E}$ such that the boundary values $\cP(z)$ are projections onto $\ran F(z)$ for almost all $z\in\T$. 

Moreover, the norm of $G$ can be estimated $\|G\|_\infty\le \delta^{-1}\|\cP\|_\infty$. 
\end{lm}
\begin{proof}
As in the proof of Nikolski's Lemma, if $G\in H^\infty_{E\shto E_*}$ is a left inverse of $F$, then $\cP:= FG$ is a bounded analytic projection. 

Le us now assume that there exists $\cP$ as in the assumptions of the lemma. Let $F^\dagger$ be a bounded left inverse of $F$ (it is always possible to find $F^\dagger$ such that $\|F^\dagger\|_\infty\le 1/\delta$). Define $G:=F^\dagger \cP$. 

It is clear that $GF=I$ a.e.~on $\T$ and that $G\in L^\infty_{E\shto E_*}$ and, moreover, 
$$
\|G\|_\infty\le \|F^\dagger\|_\infty\|\cP\|_\infty \le \delta^{-1} \|\cP\|_\infty.
$$
It remains to show that $G\in H^\infty_{E\shto E_*}$.  

Take an arbitrary $f\in H^2_E$ and let $g:= Gf = F^\dagger \cP f$. Note that  
$$
\cP f \in \cE := \{ h\in H^2_E : h(z) \in \ran F(z) \text{ a.e.~on }\T \}.  
$$
Since  $F$ is co-outer, Lemma  \ref{l-1.3} implies that $\cE =FH^2$, so $\cP f = F g$, where $g\in H^2_{E_*}$. Since $F$ is left invertible in $L^\infty_{E_*\shto E}$, such a $g$ is unique, namely $g=F^\dagger \cP f$. 

So we get that $g= G f \in H^2_{E_*}$ for all $f\in H^2_E$, which means that $G\in H^\infty_{E\shto E_*}$. 
\end{proof}

\subsection{Proof of the Generalized Nikolski's Lemma}

\begin{proof}
The implication (1)$\implies$(2) is trivial. 

\subsubsection*{(2)$\implies$(3).}
Let $F\in H^\infty_{E_*\shto E}(\Omega)$  and let $G\in H^\infty_{E\shto E_*}(\Omega)$ satisfy $FGF\equiv F$.
Define $\cP_R\in H^\infty_{E\shto E}(\Omega)$  and $\cP_K\in H^{\infty}_{E_*\shto E_*}(\Omega)$ by
\begin{eqnarray*}
\cP_R(z) & = &  F(z) G(z)\\
\cP_K(z) & = & I_{E_*}-G(z)F(z).
\end{eqnarray*}

Direct computation shows that $\cP_R^2 =\cP_R$ and $\cP_K^2=\cP_K$, so the values of $\cP_R$ and $\cP_K$ are projections.  

Clearly, $\ran \cP_R(z) \subset \ran F(z)$. On the other hand, 
$$
\ran \cP_R(z)\supset \ran (\cP_R(z) F(z)) = \ran (F(z)G(z)F(z)) =\ran F(z), 
$$
so $\ran \cP_R(z) = \ran F(z)$. 

Similarly, it is easy to see that $\ran \cP_K(z)\supset \ker F(z)$, but on the other hand
$$
F(z)\cP_K (z)= F(z)- F(z)G(z)F(z) = F(z)-F(z) = 0, 
$$
so $\ran \cP_K(z)\subset \ker F(z)$, and thus $\ran \cP_K(z) =\ker F(z)$. 

\subsubsection*{ (3)$\implies$(1)}
Let  $\cP_R(z)$ and $\cP_K(z)$ be analytic projections onto $\ran F(z)$ and $\ker F(z)$ respectively. Then for each $z\in \Omega$ the space $E$ can be decomposed into the direct sum 
\begin{equation}
\label{1.2}
E=\ran \cP_R(z) \dotplus \ker \cP_R(z), 
\end{equation}
and similarly, $E_*$ can be decomposed into 
\begin{equation}
\label{1.3}
E_* = \ker\cP_K(z) \dotplus \ran\cP_K(z).
\end{equation}
Since the functions $\cP_R$ and $\cP_K$ are bounded, then the angle between \linebreak $\ran\cP_R(z)$ and $\ker\cP_R(z)$ is uniformly bounded away from $0$. The same is true for the angle between $\ker \cP_K(z)$ and $\ran \cP_K(z)$.

In the above decompositions, the operator $F(z)$ is written as
$$
F(z)=
\left(
\begin{array}{cc}
 F_0(z) &  0   \\
 0 &  0   \\   
\end{array}
\right)\qquad \forall z\in\Omega;
$$
the condition $\ran \cP_R(z) = \ran F(z)$ means that the second row is zero and the equality $\ran \cP_K(z) = \ker F(z)$ means that the second column is zero as well. 

Note that $F_0(z) = F(z)\big| \ker\cP_{K}(z)$ if we restrict the target space to $\ran \cP_R(z)= \ran F(z)$, so $\|F_0\|_\infty\le \|F\|_\infty$.

Because the angles between the subspaces in the decompositions are uniformly bounded away from $0$, and $F$ satisfies condition \eqref{C1}, then the block $F_0(z)$ is invertible and $\nm F_0(z)^{-1}\nm\le C<\infty$ for all $z\in \Omega$. Therefore, again using the fact that the angles are bounded away from $0$, we have that the function $G$ defined by 
$$
G(z) := \left(
\begin{array}{cc}
 F_0(z)^{-1} &  0   \\
 0 &  0   \\   
\end{array}
\right)\qquad \forall z\in\Omega
$$
is bounded. It is easy to see that this function 
is  a generalized inverse for $F$. 

We want to show that $G$ is analytic, and the main problem here is that the decompositions of the domain and target spaces vary with $z$. 

The analyticity is not difficult to see in the matrix case, when both spaces $E$ and $E_*$ are finite dimensional. In a neighborhood of each point $z_0\in \Omega$ one can introduce holomorphic bases  in all four subspaces $\ran\cP_R$, $\ker\cP_R$, $\ker\cP_K$, and $\ran \cP_K$.  This gives us holomorphic bases in $E$ and $E_*$ which agree with the decompositions \eqref{1.2} and \eqref{1.3}. The matrix of $F$ in these holomorphic bases has analytic entries, so  the matrix of $F_0$ and thus the matrix of $F_0^{-1}$ have analytic entries. That means the matrix of $G$ has analytic entries, so $G$ is analytic in a neighborhood of $z_0$. Since $z_0$ is arbitrary, $G$ is analytic in $\Omega$.

Essentially, the same reasoning works in the general (infinite dimensional) case.  But, to write it formally, it is more convenient to use ``coordinate free'' notation. Namely, fix a point $z_0$ in $\Omega$. Suppose, in a neighborhood $\cU$ of $z_0$, we construct invertible bounded analytic operator-valued functions $S_1$, $S_2$,
$$
S_1(z): E^1_* \oplus E^2_* \to E_*, \qquad S_2(z): E^1\oplus E^2 \to E, \qquad \forall z\in \cU, 
$$   
such that 
\begin{alignat}{3}
\label{1.4}
S_1(z)E^1_* & = \ker \cP_K(z),& \qquad  S_1(z)E^2_* &= \ran\cP_K(z),&\qquad &\forall z\in \cU,  \\
\label{1.5}
 S_2(z) E^1 &= \ran \cP_R(z),  & S_2(z)E^2 &= \ker \cP_R(z) &\qquad &\forall z\in \cU. 
\end{alignat}
Then the function $\wt F:= S_2^{-1} F S_1$ (whose values are operators from  $E^1_* \oplus E^2_*$  to $E^1\oplus E^2$) has the following block structure:
$$
\wt F(z) = \left(
\begin{array}{cc}
 \wt{F}_0(z) &  0   \\
 0 &  0   \\   
\end{array}
\right) .
$$
Observe that the function $\wt F$ is analytic and its block structure is given with respect to constant decompositions, so the block $\wt F_0$ is analytic. Invertibility of the block $F_0$ implies the invertibility of $\wt F_0$, and  the generalized inverse  $G$ constructed above can be represented in a neighborhood of $z_0$ as 
$$
G(z)= S_2(z) \left(
\begin{array}{cc}
 \wt{F}_0(z)^{-1} &  0   \\
 0 &  0   \\   
\end{array}
\right) S_1^{-1}(z) \,. 
$$
Therefore, $G$ is analytic in a neighborhood of $z_0$, and since $z_0$ is arbitrary, $G$ is analytic in $\Omega$.

To construct $S_1$ and $S_2$ let us define 
\begin{alignat*}{2}
E^1 & := \ran \cP_R (z_0), &\qquad E^2 &:= \ker\cP_R(z_0) \\
E^1_* &:= \ker\cP_K(z_0), & E_*^2 &:= \ran\cP_K(z_0)
\end{alignat*}
and put
\begin{align*}
S_2(z)(x_1\oplus x_2) &:= \cP_R (z) x_1 + (I-\cP_R(z)) x_2 , \ x_1\oplus x_2 \in E^1\oplus E^2,\quad\forall z\in\cU\\
S_1(z)(y_1\oplus y_2) &:= (I-\cP_K(z))y_1 + \cP_K(z) y_2, \ y_1\oplus y_2 \in E^1_*\oplus E^2_*,\quad\forall z\in\cU .
\end{align*}

The functions $S_1$ and $S_2$ are bounded and analytic in a neighborhood of $z_0$. We want to show that they are invertible and satisfy \eqref{1.4} and  \eqref{1.5}. 

Consider $S_2$. Since $S_2(z_0)(x_1\oplus x_2) = x_1 + x_2$ and the angle between $E^1:=\ran \cP_R(z_0)$ and $E^2:= \ker\cP_R(z_0)$ is positive, the operator $S_2(z_0)$ is invertible. On the other hand,  
\begin{eqnarray*}
S_2(z)(x_1\oplus x_2) & = & \cP_R(z_0)x_1+(I-\cP_R(z_0))x_2\\
 &   & +(\cP_R(z_0)-\cP_R(z))x_2+(\cP_R(z)-\cP_R(z_0))x_1\\
 & = & x_1+x_2+(\cP_R(z_0)-\cP_R(z))x_2+(\cP_R(z)-\cP_R(z_0))x_1,
\end{eqnarray*}
so for $z$ close to $z_0$, the operator $S(z)$ is a small perturbation of the invertible operator $S_2(z_0)$. Therefore, $S_2$ has a bounded analytic inverse in a neighborhood of $z_0$. On the other hand, $S_2(z)E^1\subset\ran\cP_R(z)$ and $S_2(z)E^2\subset \ker \cP_R(z)$, and since $\ran\cP_R(z)$ and $\ker\cP_R(z)$ are complimentary subspaces, the invertibility of $S_2(z)$ implies \eqref{1.5}. 

The proof for $S_1$ is similar and we omit it.  
\end{proof}

\subsection{Proof of Theorem \ref{gencorona}} To prove the Generalized Operator Corona Theorem (Theorem \ref{gencorona}) we simply apply Theorem \ref{t0.2} twice; once to find the bounded analytic projection $\cP_K$, and another time to find $\cP_R$.  Then the Generalized Nikolski's Lemma (Lemma \ref{GNL}) implies that there exists a generalized inverse to $F$, thus proving Theorem \ref{gencorona}.  

Necessity in the case $\rank F(z)<\infty$ immediately follows from  Proposition \ref{p3.1} below in Section \ref{s3}. 

\section{Proof of the Main Result (Theorem \ref{t0.2})}

\subsection{Preliminaries}

\begin{lm}
Let $\Pi$ be an orthogonal projection in a Hilbert space $H$. Then  $\cP$ is a projection onto $\ran \Pi$ if and only if it can be represented as 
\begin{equation}
\label{proj1}
\cP = \Pi + \Pi V (I-\Pi), 
\end{equation}
where $V\in B(H)$. 
\end{lm}
\begin{proof}
Direct computations show that for any $\cP$ defined by \eqref{proj1}, we have 
$$
\cP^2 =\cP, \quad \text{and}\quad \ran \cP=\ran\Pi, 
$$
so $\cP$ is indeed a projection onto $\ran \Pi$. 

Now let $\cP$ be a projection onto $\ran \Pi$. Define $V:= \cP-\Pi$, so 
$$
\cP =\Pi +V. 
$$
Since $\ran\cP=\ran\Pi$, we conclude that $\ran V\subset \ran\Pi$, so $V=\Pi V$. 

On the other hand, the equality
$$
\cP x=x=\Pi x \qquad \forall x\in \ran \Pi
$$
implies that $\ran \Pi \subset \ker V$, therefore $V=V(I-\Pi)$. Then $V=\Pi V(I-\Pi)$ and \eqref{proj1} holds.
\end{proof}

One of the main tools that we will use is Green's Formula
$$
\int_{\T} u dm-u(0)=\frac{2}{\pi}\int_{\D}\Delta u\log\frac{1}{\abs{z}}dxdy,
$$
where $\Delta$ is the ``normalized'' Laplacian $\Delta := \db \p = \p \db= \frac14\left( \frac{\p^2}{\p x^2} + \frac{\p^2}{\p y^2} \right) $. If we denote by  $\mu$  the measure defined by
$$
d\mu := \frac2\pi \log\frac{1}{|z|} dxdy, 
$$ 
then Green's formula can be rewritten as 
\begin{equation}
\label{GF1}
\int_{\T} u \, dm - u(0) =  \int_{\D} \Delta u \,d\mu. 
\end{equation}

We will need some facts about orthogonal projections onto analytic families of subspaces (holomorphic vector bundles). It is a standard fact from complex differential geometry, that in the finite dimensional case, the condition $\Pi\p\Pi=0$ ($\Pi$ is a function whose values are orthogonal projections) is necessary and sufficient for the vector bundle (family of subspaces) $\ran\Pi(z)$ to be a holomorphic vector bundle. 

This fact is also true in the infinite dimensional case, if one correctly defines a holomorphic vector bundle. However, we do not know a good reference, so we will not use this fact. Instead it will follow automatically from our result. We will need the following simple lemma.  
\begin{lm}
\label{PdP}
Let $\Pi$ be a $\cC^2$ smooth (with respect to the operator norm) function of one complex variable, whose values are orthogonal projections in a Hilbert space. Assume that $\Pi\p\Pi = 0$. Then  
\begin{align*}
&(\p \Pi) (I-\Pi) =0, \qquad\p \Pi = (\p\Pi) \Pi = (I-\Pi)\p\Pi\qquad\text{and}\\
& \Delta \Pi := \p\pbar \Pi = (\p\Pi)(\p\Pi)^* - (\p\Pi)^*(\p\Pi).
\end{align*}
\end{lm}

\begin{rem}
\label{rem-PdP}
Since $\db\Pi = (\p\Pi)^*$, by taking conjugates we get the following identities for $\db\Pi$:
$$
(\db\Pi)\Pi = (I-\Pi) \db \Pi =0\textnormal{ and } \db\Pi =\Pi \db\Pi = (\db\Pi)(I-\Pi). 
$$
\end{rem}

\begin{proof}[Proof of Lemma \ref{PdP}]
The identities $\Pi=\Pi^2$ and $\Pi\p \Pi=0$ imply
$$
\p\Pi = \p \Pi^2 =(\p\Pi)\Pi+\Pi\p\Pi = (\p\Pi)\Pi,
$$
so $(\p\Pi)\Pi =\p \Pi$. 
The identity $(\p\Pi)(I-\Pi)=0$ follows immediately because 
$$
(\p\Pi)(I-\Pi) = \p\Pi - (\p\Pi)\Pi =\p\Pi -\p\Pi = 0. 
$$

The identity $(I-\Pi)\p\Pi = \p\Pi$ is an immediate corollary of the hypothesis $\Pi\p\Pi =0$:
$$
(I-\Pi)\p\Pi = \p\Pi -\Pi\p\Pi = \p\Pi. 
$$

Let us now prove the formula for $\Delta\Pi$. Taking the adjoints of both sides of the identity $\p\Pi = (\p\Pi)\Pi$, and using the fact that $(\p\Pi)^*=\overline\p\Pi$, we get 
$$
\overline\p \Pi = (\p\Pi)^* = \Pi (\p\Pi)^* = \Pi \db \Pi,
$$
so $\db \Pi = \Pi \db \Pi$. Applying $\p$ to both sides of this identity we find
\begin{eqnarray*}
\p\pbar\Pi & = & (\p\Pi)\pbar\Pi+\Pi\p\pbar\Pi\\
 & = & (\p\Pi)\pbar\Pi-(\pbar\Pi)\p\Pi+ (\pbar\Pi)\p\Pi+\Pi\p\pbar\Pi\\
 & =& (\p\Pi)\pbar\Pi-(\pbar\Pi)\p\Pi+\pbar\left(\Pi\p\Pi\right).
\end{eqnarray*}
Using the hypothesis that $\Pi\p\Pi=0$ and the fact that $\db\Pi= (\p\Pi)^*$, 
we get the final identity.
\end{proof}

\begin{rem}
\label{rem2}
We would like to mention some duality between $\Pi(z)$ and the complementary projection $\Pi\ti c(z):= I-\Pi(z)$. It is well known that the orthogonal complement of a holomorphic bundle is an anti-holomorphic bundle. 

So, it is not surprising that $\ran\Pi\ti c$ is an ``anti-holomorphic''  family of subspaces, meaning that $\Pi\ti c \db \Pi\ti c =0$ and all the identities in Lemma \ref{PdP} hold with $\Pi$ replaced by $\Pi\ti c$ and $\p$ replaced by $\db$. Another way of saying this is that the function $\Pi\ti c^\sharp(z):=\Pi\ti c(\bar z)^*$ satisfies $\Pi\ti c^\sharp \p \Pi\ti c^\sharp= 0$. 

To see this duality one should notice that $\p \Pi + \p \Pi\ti c =0$, and use the fact that $\Pi$ and $\Pi\ti c$ are complementary orthogonal projections. 
\end{rem}
\subsection{Set-Up}

First of all, let us note that the standard normal family argument allows us to assume that the function $\Pi$ is continuous, or even smooth up to the boundary, as long as we are getting uniform estimates on $\|\cP\|_\infty$. 

Namely, consider the functions $\Pi_r(z): =\Pi(rz)$, with $0<r<1$. If we can find $\cP$ for functions which are smooth up to the boundary (with a uniform estimate), then for each $r$ there exists a bounded analytic projection $\cP^r$ onto $\ran\Pi_r$ and $\|\cP^r\|_\infty\le C<\infty$. Let $\cP^{r_k}$ be a weakly convergent subsequence in $B(L^2_E(\T))$.  Letting $r_k\to 1$ we get that for any $\be_1, \be_2\in E$
$$
\La \cP^{r_k}(z)\be_1, \be_2\Ra \to \La \cP(z)\be_1, \be_2\Ra 
$$ 
uniformly on compact subsets of $\D$. Generally, a weak limit of projections does not need to be a projection. However, in our particular case, the range of the projections is fixed (for each $z\in\D$), so it follows from Lemma \ref{proj1} that the limit $\cP(z)$ is indeed a projection onto $\ran \Pi(z)$.  So, in what follows, we assume that $\Pi$ is $\cC^2$-smooth up to the boundary.

We already have a smooth projection, namely $\Pi$, and we want to add some function $V\in L^\infty_{E\shto E}$ to make it analytic and hope that we will be able to replace $V$ by $(I-\Pi) V \Pi$ without losing analyticity.    We want  $\Pi-V\in H^{\infty}_{E\shto E}$, therefore we will need the equality 
$$
\int_{\T}\La\Pi h_1,h_2\Ra dm=\int_{\T}\La Vh_1,h_2\Ra dm
$$
to hold for all $h_1\in H^{2}_E$ and $h_2\in (H^{2}_E)^{\perp}$.  Using Green's formula we get
$$
\int_{\T}\La \Pi h_1,h_2\Ra dm=\int_{\D}\p\db\La \Pi h_1,h_2\Ra d\mu=\int_{\D}\p \La (\db\Pi) h_1,h_2\Ra d\mu;
$$
where we used the harmonic extensions of $h_1$ and $h_2$ with $h_2$ being anti-analytic and $h_2(0)=0$.

By Lemma \ref{PdP}, $\db \Pi=\Pi(\db \Pi) (I-\Pi)$.  So if we define $\xi_1:=(I-\Pi)h_1$ and $\xi_2:=\Pi h_2$, then 
\begin{align*}
\int_{\D}\p\La(\db\Pi) h_1,h_2\Ra d\mu & =\int_{\D}\p\La\Pi(\db\Pi)(I-\Pi) h_1,h_2\Ra d\mu  \\
& =\int_{\D}\p\La(\db\Pi)\xi_1,\xi_2\Ra d\mu:=L(\xi_1,\xi_2).
\end{align*}
Note that the bilinear form $L$ is a Hankel form, i.e.,~$L(z\xi_1, \xi_2) = L(\xi_1, \bar z\xi_2)$. 

Suppose that we are able to prove that the estimate 
\begin{equation}
\label{main_est}
\abs{L(\xi_1,\xi_2)}\leq C\norm{\xi_1}_2\norm{\xi_2}_2,
\end{equation}
holds for all $\xi_1$ and $\xi_2$ of form $\xi_1:=(I-\Pi)h_1$, $h_1\in H^{2}_E$, $\xi_2:=\Pi h_2$, $h_2\in (H^{2}_E)^\perp$. Then, applying to the Hankel form $L$ an appropriate version of Nehari's Theorem (see Theorem \ref{t5.1} below), we get that there exists an operator-valued function $V\in L^\infty_{E\shto E}$ such that
$$
L(\xi_1,\xi_2)=\int_{\T}\La V\xi_1,\xi_2\Ra dm.
$$
Recalling the definition of $L$, $\xi_1$ and $\xi_2$ we get 
$$
\int_{\T}\La \Pi V(I-\Pi)h_1,\Pi h_2\Ra dm=L(\xi_1, \xi_2)=
\int_{\T}\La\Pi h_1,h_2\Ra dm
$$
for all $h_1\in H^2_E$ and all $h_2\in (H^2_E)^\perp$. 

Therefore, $\cP:= \Pi- \Pi V (I-\Pi)\in H^\infty_{E\shto E}$, and by Lemma \ref{proj1}, $\cP(z)$ is a projection onto $\ran \Pi(z)$ for almost all $z\in \T$.  The function $\cP$ is defined on $\T$, but we can take its harmonic extension to get a bounded analytic function in the unit disk $\D$. The identity $\cP^2=\cP$ on the boundary implies the same for all $z\in \D$, so $\cP$ is indeed a bounded analytic projection.  

To prove that $\ran\cP(z) =\ran\Pi(z)$ for all $z\in \D$, we use the inequality 
$$
\int_\D \Delta \vf\nm\xi_1\nm^2 d\mu\le C \|\xi_1\|^2_2 \qquad \forall \xi_1 := \Pi\ti c h =(I-\Pi) h, \quad h\in H^2_E, 
$$
which will be proved below in Lemma \ref{l2.3}.  The equality 
$$
\ran\cP(z) =\ran\Pi(z)\qquad \forall z\in \D
$$ 
does not immediately follow from the formula $\cP=\Pi V (I-\Pi)$.  This formula holds only on the boundary $\T$; it is completely unclear why the harmonic (analytic) extension of $\cP$ from $\T$ to the disk $\D$ is given by the same formula (with appropriately chosen  $V$).

This inequality implies that if  $h\in H^2_E$ satisfies $h(z)\in \ran\Pi(z)$ a.e.~on $\T$, then $\xi_1:= (I-\Pi)h=0$ a.e.~on $\T$, and so $\xi_1(z) = (I-\Pi(z))h(z) =0$ for all $z\in \D$, i.e.,~$h(z)\in \ran\Pi(z)$ for all $z\in \D$. That means $\ran\cP(z)\subset\ran \Pi(z)$. 

If $\dim E<\infty$, (and even if $\dim E_*<\infty$), that would be enough, because the standard reasoning with minors shows that 
$$
\dim\ran\Pi(z) = \dim\ran\cP(z)=\text{a constant}\quad \forall z\in \D \ \text{and a.e.~on }\T, 
$$ 
so the inclusion is in fact equality. But, in the general case, we need one extra step. 

Namely, $\mathcal{Q}:=I- \cP^*$ is a bounded \emph{anti-analytic} projection such that $\ran\mathcal{Q} = \ran \Pi\ti c(z)$ a.e.~on $\T$, where $\Pi\ti c :=I-\Pi$.  If we make the change of variables $z\mapsto \bar z$, we arrive at the situation considered before, so 
$$
\ran\mathcal{Q}(z)\subset \ran\Pi\ti c(z)\qquad \forall z\in \D
$$ 
(we also use the fact that $\nm\p\Pi(z)\nm = \nm\db\Pi(z)\nm = \nm\p\Pi\ti c(z)\nm = \nm\db\Pi\ti c(z)\nm$ $\forall z\in \D$). 

But 
$$
\ran\mathcal{Q}(z) = \ran(I- \cP^*(z)) = (\ran\cP(z))^\perp 
$$
{and} 
$$\ran\Pi\ti c(z) =(\ran \Pi(z))^\perp,
$$ 
so the only situation where all the inclusions are possible is $\ran \cP(z)=\ran \Pi(z)$ and $\ran\cP^*(z) = \ran \Pi\ti c (z) $ for all $z\in \D$.

\subsection{An Embedding Theorem on Holomorphic Vector Bundles}

To prove the main estimate \eqref{main_est}, we will need an  embedding theorem for functions of the form $\xi_1 =(I-\Pi) h$, $h\in H^2_E$ and $\xi_2 = \Pi h$, with $h\in \left(H^{2}_E\right)^\perp$. Such functions are not harmonic, so the Carleson Embedding Theorem does not apply. 
\begin{lm}
\label{l2.3}
Let $\vf$ be a non-negative bounded subharmonic function in $\D$ satisfying
$$
\Delta \vf(z) \geq\nm\p \Pi(z)\nm^2, \qquad \forall z\in \D, 
$$
and let $K=\|\vf\|_\infty$.  Then for all $\xi_1$ of the form $\xi_1=(I-\Pi)h$, $h\in H^{2}_E$ we have
$$
\int_\D  \Delta \vf\, \nm\xi_1\nm^2 \, d\mu \le e K e^K \|\xi_1\|_2^2 \, 
$$
and
$$
\int_\D  \nm\p \xi_1\nm^2\, d\mu \le (1 + eKe^K)\|\xi_1\|_2^2.
$$
Similarly, for all $\xi_2$ of the form $\xi_2=\Pi h$ with $\bar h\in H^{2}_E$ we have
$$
\int_\D  \Delta \vf\, \nm\xi_2\nm^2 \, d\mu \le e K e^K \|\xi_2\|_2^2 \, 
$$
and
$$
\int_\D  \nm\db \xi_2\nm^2\, d\mu \le (1 + eKe^K)\|\xi_2\|_2^2.
$$
\end{lm}

\begin{proof}
We first prove this lemma for $\xi_2$, then indicate how this proof could be used to obtain the result for $\xi_1$.  Let us take an arbitrary bounded, non-negative, subharmonic $\vf$ and  compute $\Delta \left( e^\vf \nm\xi_2\nm^2\right)$. Lemma \ref{PdP} implies that $\Pi \p \Pi =0$  and $(\p \Pi )\Pi = \p \Pi$. Therefore, using $\p h= 0$ we get  $\p \xi_2 = \p \left( \Pi h \right) = (\p \Pi) h + \Pi \p h = (\p \Pi) h = (\p \Pi) \xi_2$, and so
$$
\La \p\xi_2 , \xi_2\Ra = \La \p\xi_2 , \Pi \xi_2\Ra = \La (\p\Pi) \xi_2 , \Pi \xi_2\Ra = 0.
$$ 

Therefore,
\begin{align*}
\p \left( e^\vf \nm\xi_2\nm^2\right) & = e^\vf (\p \vf) \nm\xi_2\nm^2 + e^\vf \La \p \xi_2, \xi_2\Ra + e^\vf \La \xi_2, \db \xi_2\Ra \\
& = e^\vf (\p \vf) \nm\xi_2\nm^2 + e^\vf \La \xi_2, \db \xi_2\Ra. 
\end{align*}
Taking $\db$ of this equality (and again using $\La \xi_2 ,\p\xi_2\Ra = 0$) we get 
$$
\Delta \left( e^\vf \nm\xi_2\nm^2\right) = e^\vf \left( \Delta \,\vf \nm\xi_2\nm^2 + \nm(\db \vf) \xi_2 + \db \xi_2 \nm^2 + \La \xi_2, \Delta \xi_2 \Ra \right). 
$$
To handle $\La \xi_2, \Delta \xi_2 \Ra$ we take the $\p$ derivative of the equation $\La \xi_2 , \p \xi_2 \Ra =0$ to get 
$$
\La \p\xi_2, \p\xi_2 \Ra + \La \xi_2, \db\p \xi_2 \Ra =0, 
$$
and therefore $ \La \xi_2, \Delta \xi_2 \Ra = - \nm \p\xi_2 \nm^2 = -\nm(\p \Pi) \xi_2 \nm^2$. 
Since $\vf\ge0$
\begin{align}
\notag
& \int_{\D}\left({\Delta \varphi}\,\nm\xi_2\nm^{2}-\nm(\p\Pi) \xi_2\nm^{2}\right) d\mu \le \\
\label{2.1}   &
 \qquad\qquad \int_{\D}\left({\Delta \varphi}\nm\xi_2\nm^{2} -\nm(\p\Pi) \xi_2\nm^{2} + \nm\pbar\varphi\xi_2+\pbar\xi_2\nm^{2} \right)e^{\varphi}d\mu 
 \\   \notag
 & \qquad\qquad\qquad\qquad  = 
\int_{\T}e^{\varphi}\nm\xi_2\nm^{2}dm - e^{\vf(0)} \nm\xi_2(0)\nm^2 \le 
\int_{\T}e^{\varphi}\nm\xi_2\nm^{2}dm;
\end{align}
the  equality is just Green's formula.  In the last inequality, replacing $\vf $ by $t\vf$ with $t >1$, we get 
$$
\int_{\D}\left( t{\Delta \varphi}\nm\xi_2\nm^{2}-\nm(\p\Pi) \xi_2\nm^{2}\right)\,d\mu
\le 
\int_{\T}e^{t\varphi}\nm\xi_2\nm^{2}dm \le e^{tK} \|\xi_2 \|_2^2.
$$
Now we use the inequality $\Delta \vf \geq \nm\p \Pi\nm^2$.  It implies $\Delta \vf \,\nm\xi_2\nm^2 - \nm\p \Pi \xi_2 \nm^2 \geq 0$, and therefore
$$
(t-1) \int_{\D} {\Delta\varphi}\, \nm\xi_2\nm^{2}d\mu 
\le e^{tK} \|\xi_2\|_2^2 . 
$$
Hence,
$$
\int_{\D} {\Delta\varphi}\, \nm\xi_2\nm^{2} d\mu \le
\min_{t> 1} \frac{e^{tK}}{t-1} \,  \|\xi_2\|_2^2 = e K e^{K} \|\xi_2\|_2^2
$$
(minimum is attained at $t=1+1/K$), and thus the first statement is proved. 

To prove the second statement, put $\vf \equiv 0$ in \eqref{2.1} (we do not use any properties of $\vf$ except that $\vf \ge0$ in \eqref{2.1}) to get 
$$
\int_\D \left( \nm\db \xi_2 \nm^2 - \nm(\p \Pi) \xi_2 \nm^2 \right) \, d\mu  \le 
\int_\T \nm\xi_2\nm^2 \, dm = \|\xi_2\|_2^2. 
$$
But, the second term can be estimated as
$$
\int_\D \nm(\p \Pi) \xi_2 \nm^2 \, d\mu \le \int_\D \Delta \vf \, \nm \xi_2 \nm^2 \, d\mu \le eK e^K\|\xi_2\|_{2}^{2}
$$
and therefore $\int_\D \nm\db \xi_2 \nm^2 \, d\mu \le (1+ eK e^K) \|\xi_2\|_2^2$ .

To prove the first half of this lemma, let us  consider the complimentary projection $\Pi\ti c := I-\Pi$. The projection $\Pi\ti c$ is an orthogonal projection onto an anti-holomorphic vector bundle, see Remark \ref{rem2}. Repeating the proof for $\xi_2$ given above, but  with $\p$ and $\pbar$ interchanged,  we get the estimates for $\xi_1$. 

Alternatively, to get the estimates for $\xi_1:=\Pi\ti c h_1$, one can make the change of variables $z\mapsto \overline z$ to arrive at the situation we have treated before and use the trivial observation  
 $\nm\p\Pi(z)\nm = \nm\db\Pi(z)\nm = \nm\p\Pi\ti c(z)\nm = \nm\db\Pi\ti c(z)\nm$. 
\end{proof}

\subsection{Main Estimate \eqref{main_est}}
Now we need to estimate $L(\xi_1,\xi_2)$.  Computing the derivative inside the integral gives
\begin{eqnarray*}
L(\xi_1,\xi_2) & = & \int_{\D}\p\La (\db\Pi)\xi_1,\xi_2\Ra d\mu\\
 & = & \int_{\D}\La(\p\db\Pi)\xi_1,\xi_2\Ra d\mu+\int_{\D}\La(\db\Pi)\p\xi_1,\xi_2\Ra d\mu+\int_{\D}\La (\db\Pi)\xi_1,\db\xi_2\Ra d\mu\\
 & := & \textnormal{I}(\xi_1,\xi_2) + \textnormal{II}(\xi_1,\xi_2) + \textnormal{III}(\xi_1,\xi_2).
\end{eqnarray*}
We want to estimate each of these integrals using Lemma \ref{l2.3}.  We begin by examining the first integral more closely.  By Lemma \ref{PdP},
$$
\p\db\Pi = (\p\Pi)(\db\Pi)-(\db\Pi)(\p\Pi). 
$$
Since $(\p\Pi) (I-\Pi) =0$ and $(\db \Pi)\Pi=0$ (see  Lemma \ref{PdP} and Remark \ref{rem-PdP}), we conclude that $(\p\Pi)\xi_1 =0$, $(\db\Pi)\xi_2 =0$, so taking into account that  $(\p\Pi)^*=\db\Pi$, we get
\begin{eqnarray*}
\textnormal{I}(\xi_1,\xi_2) & := & \int_{\D}\La(\p\db\Pi)\xi_1,\xi_2\Ra d\mu\\
 & = & \int_{\D}\La(\p\Pi)(\db\Pi)\xi_1,\xi_2\Ra d\mu-\int_{\D}\La(\db\Pi)(\p\Pi)\xi_1,\xi_2\Ra d\mu\\
 & = & \int_{\D}\La(\db\Pi)\xi_1,(\db\Pi)\xi_2\Ra d\mu-\int_{\D}\La(\p\Pi)\xi_1,(\p\Pi)\xi_2\Ra d\mu=0.
\end{eqnarray*}
So we only need to estimate $\textnormal{II}(\xi_1,\xi_2)$ and $\textnormal{III}(\xi_1,\xi_2)$.  Applying the Cauchy-Schwarz inequality we get
\begin{eqnarray*}
\abs{\textnormal{II}(\xi_1,\xi_2)} & = & \left\vert\int_{\D}\La(\db\Pi)\p\xi_1,\xi_2\Ra d\mu\right\vert\\
 & \leq & \int_{\D}\left\vert\La(\db\Pi)\p \xi_1,\xi_2\Ra\right\vert d\mu\\
 & \leq & \int_{\D}\nm\db\Pi\nm\nm\p\xi_1\nm\nm\xi_2\nm d\mu\\
 & \leq & \left(\int_{\D}\nm\db\Pi\nm^2\nm\xi_2\nm^2 d\mu\right)^{1/2}\left(\int_{\D}\nm\p\xi_1\nm^2 d\mu\right)^{1/2}. 
\end{eqnarray*}
Using the fact that $\nm\db\Pi\nm^2\leq\Delta\vf$, we have
\begin{eqnarray*}
\abs{\textnormal{II}(\xi_1,\xi_2)} & \le & \left(\int_{\D}\nm\db\Pi\nm^2\nm\xi\nm^2 d\mu\right)^{1/2}\left(\int_{\D}\nm\p\xi_1\nm^2 d\mu\right)^{1/2}\\
 & \le & \left(\int_{\D}\td\vf\nm\xi_2\nm^2 d\mu\right)^{1/2}\left(\int_{\D}\nm\p\xi_1\nm^2 d\mu\right)^{1/2}.
\end{eqnarray*}
Finally,  using Lemma \ref{l2.3} we obtain
$$
\abs{\textnormal{II}(\xi_1,\xi_2)} \le  \left(eKe^K\norm{\xi_2}_2^2\right)^{1/2}\left((1+eKe^K)\norm{\xi_1}_2^2\right)^{1/2},
$$
where $K=\norm{\vf}_\infty$.  The estimate of $\textnormal{III}(\xi_1,\xi_2)$ is the same and is proved  similarly.  Combining the estimates, we have
$$
\abs{L(\xi_1,\xi_2)} \leq  2\left(eKe^K + e^2 K^2 e^{2K}\right)^{1/2}\norm{\xi_1}_2\norm{\xi_2}_2.
$$
Therefore we have that $L(\xi_1,\xi_2)$ is a bounded bilinear form.\hfill\qed

\subsection{Applying a Nehari Theorem}

We have proved that $L$ is a bounded bilinear (more precisely, a sesquilinear)  form, and we had mentioned before that $L$ is a \emph{Hankel form}, meaning that $L(z\xi_1, \xi_2 ) = L(\xi_1, \bar z \xi_2)$. Now we are going to show that one can find a symbol of this Hankel form, i.e.,~that there exists a function $V\in L^\infty_{E\shto E}$, $\|V\|_\infty\le \|L\|$ such that 
$$
L(\xi_1, \xi_2)= \int_\T \La V\xi_1, \xi_2\Ra dm 
$$
for all $\xi_1=(I-\Pi) h_1$, $h_1\in H^2_E$ and for all $\xi_2 = \Pi h_2$, $h_2\in (H^2_E)^\perp$. 

While it is possible to transform the problem so that one can apply the classical vectorial Nehari Theorem, we have at our disposal a theorem by the first author and A. Volberg, see \cite{TrlVol}, that can be applied directly to our situation. 

Let us state this theorem. Let $\mathcal{H}_1$ and $\mathcal{H}_2$ be two separable Hilbert spaces, let $S_1$ be an expanding operator ($\|S_1 x\|\ge\|x\|$) in $\mathcal{H}_1$, and $S_2$ a contractive operator ($\|S_2\|\le 1$) in $\mathcal{H}_2$ (for our problem at hand, we will have that $S_1$ and $S_2$ are isometries).  We are given an orthogonal decomposition of $\mathcal{H}_2=\mathcal{H}_2^+\oplus\mathcal{H}_2^-$, and let $S_2\mathcal{H}_2^+\subset\mathcal{H}_2^+$.  Let $\mathbb{P}_+$ and $\mathbb{P}_{-}$ be orthogonal projections in $\mathcal{H}_2$ onto $\mathcal{H}_2^+$ and $\mathcal{H}_2^-$ respectively.  Then a generalized Hankel operator $\Gamma$ is a bounded linear operator from $\mathcal{H}_1$ to $\mathcal{H}_2^-$ satisfying the following relation
\begin{equation}
\label{HankRel}
\Gamma S_1 f = \mathbb{P}_{-} S_2\Gamma f\qquad\forall f\in\mathcal{H}_1.
\end{equation}
A bounded operator $T:\mathcal{H}_1\to\mathcal{H}_2$ satisfying the commutation relation $TS_1=S_2T$ is called a generalized multiplier.  Projections of these generalized multipliers give examples of generalized Hankel operators.  The theorem of interest in \cite{TrlVol} is:
\begin{thm}(Treil, Volberg \cite{TrlVol})
\label{t5.1}
Let $S_1$ be an expanding operator and $S_2$ be a contraction.  Given a generalized Hankel operator $\Gamma$, there exists a generalized multiplier $T$ (an operator $T:\mathcal{H}_1\to\mathcal{H}_2$ satisfying $S_2T=TS_1$) such that $\Gamma=\Gamma_T$ and moreover $\norm{\Gamma}=\norm{T}$.
\end{thm}

We apply this theorem to
\begin{align*}
\mathcal{H}_1 & :=\clos_{L^2_E}\{(I-\Pi)h:h\in H^{2}_E\},\\
\mathcal{H}_2 &:=\clos_{L^2_E}\{\Pi h: h\in L^{2}_E\} \quad \text{and}\\
 \mathcal{H}_2^{-}& :=\clos_{L^2_E}\{\Pi h:h\in (H^{2}_E)^\perp\}.
\end{align*}
The operators $S_{1}$ and $S_{2}$ are defined by $S_1:=M_z\vert_{\mathcal{H}_1}$ and $S_2:=M_z\vert_{\mathcal{H}_2}$ where $M_z$ is simply multiplication by the independent variable $z$.  Then clearly, $S_2^* =M_{\bar z}\vert_{\mathcal{H}_2}$ and $S_2^* \cH_2^-\subset \cH_2^-$ so $S_2\cH_2^+\subset \cH_2^+$. 

The bilinear form $L$, defined initially on a dense subset of $\cH_1\times \cH_2^-$, gives rise to a bounded linear operator $\Gamma:\cH_1\to \cH_2^-$, $L(\xi_1, \xi_2) =(\Gamma \xi_1, \xi_2)$. We want to show that $\Gamma$ is a generalized Hankel operator, so Theorem \ref{t5.1} applies. One can see that on the dense set where $L$ is initially defined
\begin{align*}
\La \Gamma S_1\xi_1, \xi_2\Ra &= L (z\xi_1, \xi_2)  = L(\xi_1, \bar z \xi_2) \\ &= \La\Gamma\xi_1, S_2^* \xi_2\Ra 
   = \La S_2 \Gamma \xi_1, \xi_2\Ra = \La\bP_- S_2\Gamma\xi_1, \xi_2\Ra,
\end{align*}
which means that the relation \eqref{HankRel} holds, i.e.,~that $\Gamma $ is a generalized Hankel operator. By Theorem \ref{t5.1}, it can be extended to a multiplier $T: \cH_1\to \cH_2$. As one can easily see, in our case, any such multiplier is multiplication by a bounded operator-valued function $V$, whose values $V(z)$, $z\in \T$ are bounded operators from $\ran (I-\Pi(z))$ to $\ran \Pi(z)$. Of course, we can always extend $V(z)$ to the operators $V(z):E\to E$ by defining $V(z)$ to be zero on $(\ran(I-\Pi(z)))^{\perp}$.

\section{Necessity of Conditions \eqref{0.1}, \eqref{0.2}, and \eqref{0.3}}
\label{s3}
\begin{prop}
\label{p3.1}
Let $F\in H^\infty_{E_*\shto E}$ satisfy
\begin{equation}
\label{3.1}
\nm F(z)\be \nm\ge\delta \nm\be\nm, \qquad \forall z\in \D, \ \forall \be\in \left( \ker F(z) \right)^\perp
\end{equation}
for some $\delta>0$. If $\dim  \ran F(z)<\infty$, then the orthogonal projection $\Pi_R(z) $ and $\Pi\ti K(z)$ onto $\ran F(z)$ and $\ker F(z)$ respectively, satisfy conditions \eqref{0.2} and \eqref{0.3} (and therefore condition \eqref{0.1} of Theorem \ref{t0.2}). 
\end{prop}

This proposition implies Proposition \ref{p0.5}. If $\dim\ran\Pi(z)<\infty$,  Proposition \ref{p0.5} is clearly a particular case of the above proposition. If 
$$
\dim\ran\cP(z)=\infty\quad \text{and}\quad \codim\ran\cP(z)<\infty,
$$ 
we should consider  the anti-analytic projection $I-\cP(z)^*$. Note that $\ran(I-\cP(z)^*) = \ran\cP(z)^\perp$ is finite dimensional, so
applying the ``anti-holomorphic version'' of Proposition \ref{p3.1} will yield the result. 

Namely,  $\ran \Pi\ti c (z)$ is an anti-holomorphic vector bundle, see Remark \ref{rem2},  so applying to it an ``anti-holomorphic version'' of Proposition \ref{p3.1} (or applying the proposition to the holomorphic projection $I-\cP(\overline z)^*$), we get that the complimentary orthogonal projection $\Pi\ci c := I-\Pi$ satisfies conditions \eqref{0.2} and \eqref{0.3}, and thus condition \eqref{0.1}. But
$\p \Pi = -\p \Pi\ti c$, so $\Pi$ satisfies these conditions as well.

\begin{proof}[Proof of Proposition \ref{p3.1}]
Let us first consider the case when $\ker F(z)=\{0\}$ (for all $z\in\D$) and  $\dim E_*<\infty$. Assumption \eqref{3.1} implies that $F^*F$ is invertible in $L^\infty(\D)$, so one can write the formula for the orthogonal projection $\Pi(z)$ onto $\ran F(z)$, 
$$
\Pi := F(F^*F)^{-1}F^*.
$$
Direct computations show that 
$$
\p \Pi = (I-\Pi) F' (F^*F)^{-1}F^*, 
$$
so 
\begin{equation}
\label{3.2}
\nm \p \Pi(z)\nm \le C \nm F'(z)\nm.
\end{equation}
It is well known and easy to show that for a scalar function $F\in H^\infty$,
$|F'(z)|\le C(1-|z|)^{-1}$ and the measure $|F'(z)|^2 (1-|z|) dx dy$ is Carleson. Similarly, for a bounded analytic function in the disk $\D$ taking values in a Hilbert space%
\begin{align}
	&\nm F'(z)\nm\le C/(1-|z|); \hfill\hfill\hfill 
	\label{3.3}\\
  &\text{the measure } \nm F'(z) \nm^2 (1-|z|) dxdy \text{ is a Carleson measure. \hfill}     \label{3.4}
\end{align}
These estimates are not true for functions with values in arbitrary Banach spaces. It is very easy to construct a counterexample for functions with values in $\ell^\infty$.
%
%

Since $\dim E_*<\infty$ our operator-valued function $F$ is a bounded analytic function taking values in the Hilbert-Schmidt class $\mathfrak{S}_2$, $\nm F(z)\nm_{\mathfrak{S}_2} \le (\dim E_*)^{1/2}\nm F(z)\nm$. But $\mathfrak{S}_2$ is a Hilbert space, therefore \eqref{3.3} and \eqref{3.4} and thus \eqref{0.2} and \eqref{0.3} hold.

Condition \eqref{0.1} follows from \eqref{0.2} and \eqref{0.3}, but also can be verified directly: one can put $\vf:= C\tr (F^*F)$. This is probably the easiest possible $\vf$ to find, but it does not give a good estimate in the matrix Corona Theorem. The function $\vf = \ln\det(F^*F)$ gives the best known estimate, see \cite{TreilWick}. 

To treat the general case, let us first note that the assumptions of the propositions imply that $\dim \ran F(z) \equiv \text{Const}<\infty$ for all $z\in \D$. Indeed, \eqref{3.1} and the continuity of $F$ imply that $\ran F(z) $ depends continuously on $z$, so its dimension is constant.  

Following ideas of M.~Andersson, see \cite{Andersson1}, consider the orthogonal projection $P(z)=\Pi\ci K(z)$  onto $\ker F(z)$.  The function  $F^*F +P^*P$ is invertible in $L^\infty(\D)$ (note that $P^*P=P^2 =P$, but we write $P^*P$ to keep the notation symmetric), and 
$$
\Pi := F(F^*F + P^*P)^{-1}F^*. 
$$
Direct computations (using the fact that $FP=0$) show that 
\begin{equation}
\label{3.5}
\p \Pi = (I-\Pi) F'(F^*F + P^*P)^{-1} F^*
\end{equation}
so \eqref{3.2} holds. Repeating the  reasoning for the case $\ker F(z) \equiv \{0\}$, we get that conditions \eqref{0.1}, \eqref{0.2}, and \eqref{0.3} hold. 

There is a small detail of note here; to prove \eqref{3.5}, one needs to assume that $P$ is differentiable.  However, it is easy to take care of this problem. Namely, let  $F^*(z_0)\be_k$, $\be_k\in E_*$, $k=1, 2, \ldots, d$ be a basis in $\ran F^*(z_0)$ ($\dim\ran F(z)=\dim \ran F^*(z)$). Then for $z$ close to $z_0$, the vectors $F^*(z)\be_k$, $k=1, 2, \ldots , d$ form a basis in $\ran F^*(z)$. If we let $R$ denote the restriction of $F^*$ to $\spn\{\be_k: k=1, 2, \ldots , d\}$, then in a small neighborhood of $z_0$ $\ker R(z) =\{0\}$ and $\ran R(z)= \ran F^*(z)$. Since $\ker F(z) =(\ran F^*(z))^\perp$, we get the formula for the orthogonal projection $P(z)$ onto $\ker F(z)$, 
$$
P:= I- R(R^*R)^{-1} R^*, 
$$
so $P$ is clearly differentiable. 

To prove the statement about $\Pi\ti K$, we note that $\Pi\ti K$ can be represented as 
$$
\Pi\ti K := I - F^* (FF^* + \Pi\Pi^*)^{-1}F
$$
and repeat the reasoning for $\Pi$. Notice, that differentiability of $\Pi$ is already proved,  so we do not need to worry about it. 
\end{proof}


\def\cprime{$'$}
\providecommand{\bysame}{\leavevmode\hbox to3em{\hrulefill}\thinspace}

\end{document}